\documentclass[12pt,letterpaper]{amsart}    

\usepackage[parfill]{parskip}    % Activate to begin paragraphs with an empty line rather than an indent

\usepackage{graphicx}

\usepackage{amssymb} %math symbols

\usepackage{amsmath} %math stuff

\usepackage{amsthm} %Theorems and Such

\usepackage{mathrsfs}

\usepackage{enumerate} %allows to change bullet points in enumerate

\usepackage{fullpage}  

\usepackage{setspace} %Spacing

\usepackage{color}

\usepackage{wrapfig}

\usepackage{url}

\theoremstyle{plain}
\begin{document}

\title{Distinguished graduates in mathematics of Jagiellonian University in the interwar period. Part I:  1918-1925}

\author[S. Domoradzki]{Stanis\l aw Domoradzki}
\address{Faculty of Mathematics and Natural Sciences,  University of Rzesz\'ow,  ul. Prof. S. Pigonia 1, 35-959 ,  Rzesz\'ow, Poland}
\email{domoradz@ur.edu.pl}
\author[M. Stawiska]{Ma{\l}gorzata Stawiska}
\address{Mathematical Reviews, 416 Fourth St., Ann Arbor, MI 48103, USA}
\email{stawiska@umich.edu}

\date{\today}                                           % Activate to display a given date or no date

\maketitle

%\newpage

 \tableofcontents

\setcounter{tocdepth}{3}

\begin{abstract}
In this study, we present profiles of some distinguished graduates in mathematics of the Jagiellonian University from the years 1918-1939. We discuss their professional paths and scholarly achievements, instances of scientific collaboration, connections with other academic centers in Poland and worldwide, involvement in mathematical education and teacher training,  as well as their later roles in
Polish scientific and academic life.  We also try to understand in what way they were shaped by their studies and how much of Krak\'ow scientific traditions they continued.  We find  strong support for the claim that there was a distinct, diverse  and deep mathematical stream in Krak\'ow between the wars, rooted in classical disciplines like differential equations and geometry, but also open  to  new trends in mathematics. 
Part I concerns the graduates before the university reform, in 1918-1926. \\
%Part II concerns people who graduated after the university reform, getting the newly introduced master's degree %in mathematics, 1926-39.
\par W niniejszej pracy przedstawiamy sylwetki niekt\'orych wybitnych absolwent\'ow Uniwersytetu Jagiello\'nskiego w zakresie matematyki z lat 1918-1939. Omawiamy ich drogi zawodowe i osi\c agni\c ecia naukowe, przyk\l ady wsp\'o\l pracy naukowej, zwi\c azki z innymi o\'srodkami akademickimi w Polsce i na \'swiecie, zaanga\.zowanie w nauczanie matematyki i kszta\l cenie nauczycieli oraz ich p\'o\'zniejsze role w polskim \.zyciu  akademickim. Pr\'obujemy tak\.ze zrozumie\'c, w jaki spos\'ob zostali oni ukszta\l towani przez swoje studia  i na ile  kontynuowali krakowskie tradycje naukowe. Znajdujemy mocne dowody na poparcie tezy, \.ze w Krakowie mi\c edzywojennym  istnia\l \ wyra\'zny, zr\'o\.znicowany i g\l \c eboki nurt matematyczny, zakorzeniony w dyscyplinach klasycznych takich jak r\'ownania r\'o\.zniczkowe i geometria, ale r\'ownie\.z otwarty na nowe trendy w matematyce. Cz\c e\'s\'c I dotyczy  absolwent\'ow sprzed reformy kszta\l cenia uniwersyteckiego, za lata 1918-1926.

\end{abstract}

\section{Introduction}

Mathematical traditions at Jagiellonian University in Krak\'ow (Cracow) are several  centuries old. Founded in 1364, the university got its first chair for mathematics and astronomy in 1405. It functioned during the rule of the Jagiellonian dynasty and elected kings, as well as during the 123-year period  of occupation by the Austro-Hungarian empire  (\cite{Go64}, \cite{Dy00}, \cite{CP12}).  World War I brought some  interruptions to its activity.  Krak\'ow was a fortress,  important for the war operations. In anticipation of a Russian attack, emergency evacuation was ordered and one-third of the  population left the city.  A few university  buildings, including its clinics,   were requisitioned for the military needs. The university was closed for the winter semester of the academic year 1914/15. However, the the offensive was thwarted before the Russians approached  Krak\'ow  and soon  the classes resumed. They  continued until the bloodless liberation of Krak\'ow on October 31, 1918, when Polish soldiers disarmed an Austrian garrison stationed in the City Hall Tower and the occupying army capitulated. To commemorate  the event, the words ``Finis Austriae" were entered in the book of doctoral promotions at the Jagiellonian University.  (\cite{Do15}) \\
On November 11, 1918 (the same day the armistice  was signed in Compi\`egne), the Regency Council handed in the military command to J\'ozef Pi\l sudski. On November 14, it passed all remaining authority to him and dissolved itself. On November 16, Pi\l sudski notified the Entente Powers about the emergence of the independent Polish state.\\
The beginnings of the newly proclaimed Republic of Poland were turbulent. Even though the Great War was over, the state had to fight for its borders against the competing interests of Germans as well as Czechs, Slovaks,  Lithuanians and Ukrainians, who also aspired to independence, and finally, in 1919-1921, to defend itself against the Bolshevik Russia.  Many students  interrupted their courses to volunteer for the military service. Internally, the emerging republic also faced many difficulties. It had to forge, among other things, a common legal, administrative, monetary  and educational systems, in place of no longer relevant  disparate systems of the three occupying powers. The existing Polish-language academic institutions-- that is, the Jagiellonian University,  the Lw\'ow (Lvov) University and the Lw\'ow Polytechnics-- played a major role in these endeavors, educating the future forces to staff Polish administration, courts and schools at all levels, including the new or revived academic schools in Polish territories formerly under Russian or German occupation. \footnote{Since 1871, Polish was  the language of instruction in the institutions of  higher education in the Polish province of Galicia of the Austro-Hungarian empire.}  Of course hardships of everyday life affected the educational process and the scientific research. At the opening of the academic year 1922/23, the vice-rector of the Jagiellonian University, Professor Stanis\l aw Estreicher, a historian of law and a bibliographer, characterized the current conditions as follows:\\

\textit{``It cannot be denied (...) that the current relations, resulting from  the tremendous shock of civilization which was the Great War, do not create by any means conditions that are favorable for the development of scientific institutions. The funding of our institutes and laboratories, which was not great even before the war,  in current conditions does not bear any proportion even to the most modest needs.  Every year we lose more and more scientific contact with abroad, we are more and more isolated from intellectual trends, foreign books and  journals reach our hands less and less frequently, with more and more difficulty. Even the most important scientific institutions  lack  the space or are unable to expand, and   it is almost impossible to bring scientific staff from outside because of the difficulties in finding accommodations for them. The same difficulty affects  the young people: our youth struggle not only with provisional difficulties and  with dearness of scientific aids, but primarily with the difficulty of finding accommodations in Krak\'ow. The number of youth does not increase at the rate at which it should, if one takes into account how big is the need of the society for intellectual forces, lacking in Poland in every way."} (\cite{KrUJ2122}; translated by the second author)\\

A committee  was formed to help the students. Institutions and more affluent private citizens answered the appeal and  contributed cash and commodities. Despite all difficulties, young people were trying to pursue academic   studies to satisfy their intellectual needs and  prepare themselves for their roles in the public life of the new state. In the academic year 1921/22, the total number of students at the Jagiellonian University was 4580. Most of them, 2002, were enrolled in the faculty of philosophy (where mathematics also belonged). There were 976 women students, 763 of them in the faculty of philosophy (\cite{KrUJ2122}).\\

In 1918 there were three professors of mathematics at the Jagiellonian University: Stanis\l aw Zaremba, Kazimierz \.Zorawski and Jan Sleszy\'nski. Lectures were also given by Antoni Hoborski, Alfred Rosenblatt and W\l odzimierz Sto\.zek, soon joined by  Leon Chwistek, Witold Wilkosz  and Franciszek Leja. The standard curriculum,  established before 1918, comprised mathematical analysis (differential and integral calculus), analytic geometry, introduction to higher mathematics, differential equations, differential geometry, theory of analytic functions. Less frequently,   lectures were given in number theory, power series, synthetic geometry, algebraic equations, higher algebra and projective geometry. Occasionally, there were also lectures in theory of elliptic functions, theory of transformations, theory of conformal mappings, spherical trigonometry, calculus of variations, analytic mechanics, algebraic geometry, elementary geometry, kinematics of continuous media, equations of physics, probability, theory of determinants, introduction to the methodology of mathematics,  mathematical logic, integral equations. (see \cite{Go64}) New subjects were added by Witold Wilkosz, who became an extraordinary professor in 1922: foundations of mathematics, set theory, theory of quadratic forms, theory of functions of a real variable, geometric topology, group theory. The course of studies could be concluded with either the teacher's examination or PhD examination.  An additional subject  was required at the examination (most commonly, physics was chosen). \\

A common view is that Krak\'ow was lagging behind the mathematical centers of Lw\'ow and Warsaw, neglecting  developments in topology, set theory and functional analysis, in which the Polish school of mathematics achieved its most remarkable and lasting results.  It has been claimed that the emphasis on classical mathematics, mostly on  differential equations, caused the exodus from Krak\'ow of younger mathematicians more interested in pursuing modern fields and eventually led to the dismantling of Krak\'ow research group. Here are the impressions of  the outstanding Russian mathematician  Nicolai Lusin after his visit to Warsaw, expressed in his letter to Arnaud Denjoy, dated September 30, 1926  \\

\textit{``It seems to me that the mathematical life in Poland follows two completely different ways: one of them is inclined to the classical parts of mathematics, and the other to the theory of  sets (functions). These ways exclude one another in Poland, being the irreconcillable enemies, and now a fierce struggle is going on between them."} (cited after \cite{Dem14}; translated by S. S.  Demidov)\\

\textit{``The classical side is currently represented only by the ancient (over 500-year old) Krak\'ow University and Krak\'ow Academy. Among Polish mathematicians the most stalwart supporter of this way is  Professor Zaremba. Other supporters of this tendency stick close to Mr. Zaremba. [...] Thus currently only Krak\'ow is a stronghold of classical mathematics. However, the Polish mathematicians whom I saw in Warsaw claim univocally that Mr. Zaremba's cause is doomed to failure, and this is why numerous colleagues and students are leaving him. Thus Mr. Zaremba's students, Dr. Kaczmarz and Dr. Nikliborc, are leaving Krak\'ow for Lw\'ow; Mr. Banach and Mr. Sto\.zek already did so. A student of Mr. Zaremba, Mr. Leja, left for Warsaw.  Mr. Nikodym also intends to leave Krak\'ow and only material difficulties do not allow him to relocate to a different city. One can therefore speak about  dissolution of the group of Krak\'ow mathematicians."} (cited after \cite{Dy00}; translated by the second author). \\

The truth is more complex than this. We already mentioned new subjects introduced (mainly by Wilkosz) to the students' curriculum which reflected the latest mathematical developments.  These subjects also made their way into research of Krak\'ow mathematicians, e.g., Wa\.zewski did some significant work in topology and Go\l \c ab in geometry of metric spaces (see the subsequent sections devoted to  individual  mathematicians for more detail). Let us  now address the relocations of Krak\'ow mathematicians.  In the first years of the Republic of Poland, it was necessary to provide Polish-speaking staff to the revived or newly created institutions of higher education. Lw\'ow and Krak\'ow mathematicians took the opportunity to engage in building academic and scientific life where there was none.  Kazimierz  \.Zorawski, who  worked mainly in differential geometry and  in the academic year 1917/18 served as the rector of the Jagiellonian University,  moved to Warsaw in  1919, to become a professor first at  the  Warsaw Polytechnic, then at the University of Warsaw (both reopened in 1915) and to serve as an official in the Ministry for Religious Denominations and Public Enlightenment. He remained in Warsaw until his death in 1953. Franciszek Leja, an assistant in Krak\'ow, whose works concerned functions of a complex variable, potential theory  and approximation theory,\footnote{Leja was also  first to define explicitly a topological group.} obtained a chair of mathematics at the Warsaw Polytechnics in 1923. However, he returned to Krak\'ow  to become an ordinary professor of mathematics in 1936. A few mathematicians who were based in Lw\'ow before World War I-- Wac\l aw Sierpi\'nski, Zygmunt Janiszewski and Stefan Mazurkiewicz--decided to take positions at the University of Warsaw. Their departures (as well as  the death of J\'ozef Puzyna in 1919) left  the Lw\'ow mathematics understaffed. Some of the gaps were filled by Hugo Steinhaus, W\l odzimierz Sto\.zek, Stefan Kaczmarz, W\l adys\l aw Nikliborc and Stefan Banach. Without formal training in mathematics and university diploma, Banach had no chance for a standard academic employment. The invitation from Lw\'ow to work as Antoni \L omnicki's assistant (issued at the initiative of Hugo Steinhaus) came as a unique opportunity. Steinhaus himself  did not have an academic appointment for a few years following his PhD, even though he got his habilitation in Lwów  in 1917.  After World War I he worked as a mathematical expert for a gas company and considered himself a ”private scholar” before joining the Lwów University. Sto\.zek became a professor of mathematics at the Lw\'ow Polytechnics, but soon he ceased to do research and devoted himself to writing textbooks. Kaczmarz and Nikliborc went on to get their PhD degrees from Lw\'ow and to build their mathematical careers, ended by their premature deaths (Kaczmarz fell in the September campaign of 1939 and Nikliborc was driven to suicide by Communist security forces in 1948). While Kaczmarz's research concerned real analysis, orthogonal series and numerical methods, Nikliborc continued to work on potential theory, differential equations and mechanics. Otto Nikodym held on to his (secure and relatively well-paid) job as a high-school teacher, but  ultimately left Krak\'ow to  pursue an academic career, getting his PhD (in 1924) and  habilitation (in 1927) at the University of Warsaw. One can thus say that the relocations were motivated not so much by individuals’ wishes to  pursue particular mathematical  interests, but rather by available career opportunities.\footnote{There  were also mathematicians who got their PhD in Krak\'ow and were employed by other centers, e.g.,  Stefan Kempisty (PhD 1919), later a professor at the Stefan Batory University in Wilno, working mainly in real analysis, and Witold Pogorzelski (PhD 1919), later a professor of Warsaw Polytechnics, working mainly  in differential and integral equations.} While Zaremba remained faithful to classical mathematical disciplines, he did not prevent his  colleagues from pursuing topics in real analysis, set theory, topology and other modern developments, or from teaching them to students. His main neglect was in not taking care to create positions for junior mathematicians. Here is an account by Andrzej Turowicz, a student in the years 1922-26:\\

\textit{``When I came for my first year of studies at the university, there were no assistants. So-called ``proseminarium" [a pre-seminar], which preceded the recitations, was led by docent Leja. I was in this pre-seminar of his. He was a high school professor and had contract classes [at the university]. Only when I entered the second year,  Wilkosz took care to have two deputy assistants nominated, that is Le\'sniak and Mrs. Wilkosz ([who was] a mathematician). Then also Turski made it. Even later, there was Krystyn Zaremba."}\\

\textit{``I taught high school for 10 years. At that time, there was one adiunkt and one deputy assistant at the university. That was it. There was no chance of getting into university [as junior faculty]. It was Zaremba's fault. He did not take care to acquire an assistant."} (\cite{TuAU}. \\

 In the initial years, however, the Jagiellonian University and the Academy of Mining (created in 1919) were able to accommodate the mathematical talents of Tadeusz Wa\.zewski and Stanis\l aw Go\l \c ab. They became world-class researchers and later created their own scientific schools. There were a few other distinguished graduates of mathematics in the period 1918-1926 (before the unification of university system and introduction of the master's degree for all disciplines) who later made an impact on Polish intellectual and academic life. Below we present their profiles.

\section{Profiles}

\subsection{Tadeusz Wa\.zewski (1896-1972)} 

\begin{wrapfigure}{L}{0.45\textwidth}

\begin{center}
\includegraphics[width=0.45\textwidth]{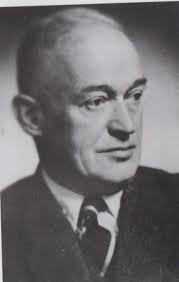}
\end{center}

\end{wrapfigure}

Born in the village of Wygnanka, in the Austro-Hungarian province of Galicia (later in Tarnopol Voivodeship in the II Republic of Poland, then in USSR, now Ukraine), he passed his \emph{matura} in the 1st Gymnasium in Tarn\'ow. In the years 1914-1920 he studied   in the philosophical faculty of the Jagiellonian University  He started studying physics, then, influenced by S. Zaremba,  switched to mathematics. Military service interrupted his studies. When Poland regained independence in November 1918, he served on a citizen patrol  in Krak\'ow.  In the years 1920-21 he taught in St. Anne State Gymnasium in Krak\'ow as a substitute teacher. In the years 1921-23 he studied at the University of Paris, where in 1924 he obtained PhD degree on the basis of the thesis ``Sur les courbes de Jordan ne renfremant aucune courbe simple ferme\'e de Jordan" concerning dendrites, i.e., locally closed continua which do not contain simple closed curves (later published as \cite{Wa23b}). On the doctoral committee there were \'Emile Borel, Arnaud Denjoy and Paul Montel (\cite{Do13}). Paul \'Emile Appell is also signed on Wa\.zewski's diploma, as the rector of the university of Paris at that time. Wa\.zewski returned to Krak\'ow, where he got his habilitation in 1927 for a thesis concerning rectifiable continua (\cite{Wa27}). He became a deputy professor in  1926 and an extraordinary professor of the Jagiellonian University in 1933 (\cite{PANWaz}, \cite{UJWaz}).\footnote{The Council of the Philosophical Faculty of Jagiellonian University, of which mathematics was part, proposed Wa\.zewski  as a candidate for extraordinary professorship already in 1929. However, in a poll among mathematicians from other centers, Wa\.zewski's name was not mentioned. Instead, other names appeared: of W\l adys\l aw Nikliborc, Otto Nikodym,  Antoni Zygmund, Bronis\l aw Knaster, Alfred Rosenblatt,  Alfred Tarski, Juliusz Schauder, Witold Hurewicz and Aleksander Rajchman.   The ministry   concluded that Wa\.zewski's candidacy was inferior and rejected the proposal (\cite{Dy00}.} \\

Wa\.zewski  was arrested along with other professors in the Sonderaktion Krakau on November 6,  1939, imprisoned in the Sachsenhausen concentration camp and released in February 1940. He stayed in Krak\'ow until the end of the Nazi occupation, officially lecturing at the Men School of Commerce while also engaging in clandestine teaching and conducting an illegal mathematical seminar with many participants who presented their results. In 1945 he became an ordinary professor and took active efforts to rebuild academic life in Krak\'ow, not only at the Jagiellonian University, but also at the State  Pedagogical College and Academy of Commerce as well as at the State Mathematical Institute (later Mathematical Institute of the Polish Academy of Sciences). He lectured, published papers,  edited journals (e.g. \emph{Annales Polonici Mathematici}), headed academic units (a university chair and a division of the Mathematical Institute), presided the Polish Mathematical Society (1959-61) and, perhaps most importantly, supervised PhD students, many of whom later became  distinguished mathematicians.  He gave   short  communications at the International Congresses of Mathematicians, in Oslo (1936) and  Stockholm (1962), and a plenary address in Amsterdam (1954)(\cite{ICM}) \\

Initially Wa\.zewski was interested in point-set topology and obtained important results in this area.  In his PhD thesis he constructed a universal dendrite, i.e., a dendrite containing a homeomorphic image of any other dendrite. This set, known as Wa\.zewski's dendrite, still attracts interest of mathematicians, either in its own right (\cite{BCMMS} or as a tool in the study of Berkovich analytic spaces, lying at the intersection of algebraic geometry, number theory and topology (\cite{HLP14}).   In  1923 Wa\.zewski proved a theorem on properties of the hyperspace of a locally connected metric continuum (\cite{Wa23a}). This result was  obtained independently and published amost simultaneously by Leopold Vietoris (\cite{Vi23}) and later generalized by Menachem Wojdys\l awski (\cite{Wo39}) and other mathematicians (\cite{Ch98}).  In later years he occasionally returned to purely topological topics, e.g. in \cite{Wa61}. However,  differential equations  became his main interest, the area of his highest achievements and a springboard for new directions of research.  He started working on them in 1930s and published over 25 papers concerning them before World War II. Two of these works he wrote jointly with Stanis\l aw Krystyn Zaremba (\cite{WaZa35}, \cite{WaZa36}), and in at least one he addressed questions inspired by the work of Adam Bielecki (\cite{Wa31}).  He also worked on differential inequalities. In 1947 Wa\.zewski published the first version of his fundamental result, later called Wa\.zewski's retract principle (\cite{Wa47})-- a creative application of topology to differential equations. The result has many variants now (\cite{Srz11}). Roughly speaking, it says that certain properties of the solutions of a differential equation (or a system of equations) on the boundary of a given domain imply that some solutions to the equation must stay in the domain. The case of systems of equations requires significant  use of the topological notion of  retract. Versions for dynamical systems were also developed and led to the theory of so-called Conley index (\cite{Ry12}), which is a tool for analyzing topological structure of invariant sets of smooth maps or smooth flows. In 1948 (\cite{Wa48}), investigating the domain of existence of an implicit function, Wa\.zewski introduced and studied a linear differential equation, later called Wa\.zewski equation. The existence of polynomial solutions to this equation was later shown (\cite{MO87}) to be equivalent to the famous Keller Jacobian Conjecture (which is still open in the complex case). In 1960s Wa\.zewski wrote a series of papers in the optimal control theory, in which he created a new direction in the area by using the notion of a differential inclusion, introduced by S. K. Zaremba and A. Marchaud (\cite{Gor11}, \cite{Ol76}). \\

Wa\.zewski had a great pedagogical talent. He was ranked very highly by students as an instructor. He published lecture notes (\cite{Wa48l}) and  in his research he sometimes took up topics of educational value. For example, he gave a unified proof of all cases of the de l'H\^ospital rule in calculus (\cite{Wa49}; later he also published a version for Banach spaces). He also had a strong sense of mission.  During the occupation he  entered twice the Krak\'ow ghetto to talk to an amateur mathematician named Rappaport (first name unknown). The meetings were arranged by Tadeusz Pankiewicz (\cite{Pa03}, \cite{Pe00b}), a Pole who was allowed by the Nazis to run a pharmacy inside the ghetto until its end. Wa\.zewski risked his own life, as he was a former concentration camp prisoner and the entrance was illegal, but he had discussions with Rappaport, who gave an approximate solution to the  angle trisection problem. Rappaport did not survive the war; his result was  published by Wa\.zewski  in 1945.  (\cite{Wa45}).  Another memorable instance of Wa\.zewski' s doing what he considered the right thing was the following: his former student Andrzej  Turowicz became a Benedictine monk and a priest (Fr. Bernard) after the war. It became hard for Turowicz to participate actively in the scientific life, as the communist authorities kept  limiting the influence of the Catholic church in the public life. But thanks to Wa\.zewski's firm stand, he was employed in the Mathematical Institute of the Polish Academy of Sciences since 1961, had his habilitation approved in  1963 and became an extraordinary professor in 1969.  (For more information about Turowicz, see \cite{DS15} and the references therein). \\

\subsection{W\l adys\l aw Nikliborc (1899-1948)} \begin{wrapfigure}{L}{0.45\textwidth}

\begin{center}
\includegraphics[width=0.45\textwidth]{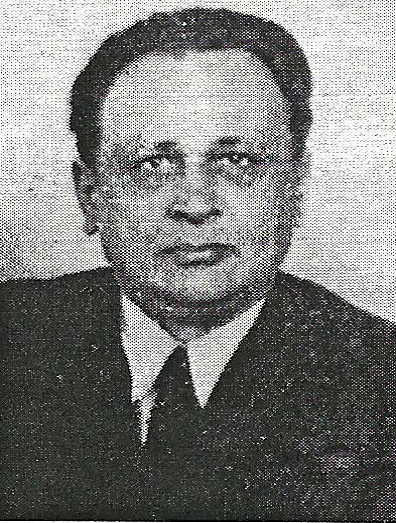}
\end{center}

\end{wrapfigure} Born in Wadowice, he finished high school there and passed his \textit{matura} examination in 1916. His university studies  were interrupted by military service. Released from the military in December 1920, he studied mathematics at the Jagiellonian University, graduating in 1922. He attended lectures in mathematics (by Stanis\l aw Zaremba, Kazimierz \.Zorawski, Jan Sleszy\'nski, Antoni Hoborski, Alfred Rosenblatt, W\l odzimierz Sto\.zek, Witold Wilkosz), physics (by W\l adys\l aw Natanson, Czes\l aw Bia\l obrzeski, Konstanty Zakrzewski, Stanis\l aw Loria) and astronomy (by Tadeusz Banachiewicz). In October 1922 he became an assistant in the Chair of Mathematics headed by Antoni \L omnicki  at the Faculty of Mechanics of  the Lw\'ow Polytechnics. He also taught mathematics at a private Ursuline gymnasium for women.  \\

In 1924 Nikliborc obtained his PhD degree  at the Faculty of Philosophy of the Lw\'ow University, passing exams in mathematics, astronomy and philosophy and presenting a thesis ``On applications of the fundamental theorem of Cauchy on the existence of solutions of ordinary differential equations to boundary value problems for the equation $y''=f(x,y,y')$."  His habilitation at Lw\'ow University took place in 1927 and was based on a two-part paper ``Sur les fonctions hyperharmoniques", published in Comptes Rendus of the Paris Academy of Sciences in 1925 and 1926 (\cite{Ni25}, \cite{Ni26}). In this paper Nikliborc considered the Dirichlet problem in a polydisk for a function which is the real part of a  holomorphic function of two complex variables. Nowadays such functions are called pluriharmonic. The study of these functions was initiated in 1899 by Henri Poincar\'e and continued by Luigi Amoroso (\cite{Am12}), who in 1912 also considered a Dirichlet problem (in a somewhat more general domain), but treated it as a problem in four real variables.  It was Nikliborc who made subsequent  advances in the theory of pluriharmonic functions. The referee of his habilitation, Hugo Steinhaus, praised not only his ingenuous methods but also his thorough knowledge of differential equations. In 1931 Nikliborc got  habilitation  in theoretical mechanics at the Lw\'ow Polytechnics, on the basis of a paper about rotating fluid. He became interested in   fluid mechanics and  celestial mechanics and  worked on these subjects. From 1937 to 1939 he was an extraordinary professor at the Warsaw Polytechnics.  He wrote academic and high school textbooks (some of them jointly with Steinhaus or W\l odzimierz Sto\.zek).  He spent the years of World War II in Lw\'ow. During the Soviet occupation (1939-1941) he was a professor of mathematics in Lw\'ow, and during the German occupation (1942-1944) he lectured at the Staatliche Technische Fachkurse, a school for vocational training created as partial replacement of the Polytechnics (which, like all Polish institutions of higher education, was closed by the Nazis). He took care of ailing Stefan Banach in his last days. In 1945  he went to Krak\'ow, to  take the Chair of Mathematics in  the Faculty of Engineering at the Academy of Mining and Metalurgy and to give lectures at the Jagiellonian University, but soon he  moved to Warsaw (Andrzej Turowicz was appointed to lecture in his stead, \cite{TuAU}).  He became an ordinary professor there: first at the Polytechnics, then at the University. He committed a suicide in 1948 in Warsaw after being arrested and interrogated by the communist security forces (\cite{Du07}).  His advanced textbook in differential equations was published posthumously in 1951.  (See also \cite{Sl48}.)\\

 Nikliborc mastered differential equations as a student in Krak\'ow. Rumor had it (as related later by Fr. Turowicz, \cite{TuAU}) that he and Stefan Kaczmarz, an earlier Jagiellonian student, moved to Lw\'ow for their PhD in order to avoid taking further examinations with Zaremba. Because of lack of materials, we cannot confirm or deny personal factors playing a role, but  both graduates were certainly well prepared for scholarly work in mathematics. While Kaczmarz's interests were closer to  those of the Lw\'ow mathematical school, Nikliborc continued working on classical analysis, differential equations and generalization of harmonic functions. However, he wrote a joint paper with Kaczmarz, \cite{KaNi28}, and entered several problems in the Scottish Book, one as a ``Theorem" (\# 128, 129, 149 and 150, \cite{KsSz3541}). His later  interests in mechanics were influenced mainly by Leon Lichtenstein, whom he visited in Leipzig in 1930-1931.  Solid background in physics and astronomy  acquired in Krak\'ow  helped Nikliborc in his work on equilibrium figures in hydrodynamics and on the three-body problem.\\

\subsection{Stanis\l aw Bilski (1893-1934(?))}

Stanis\l aw Bilski was born in Zgierz. He passed his \emph{matura} examination in 1911 in \L \'od\'z Merchants' School of Commerce. In the years 1914-1918 he worked as a teacher in Zgierz. He started his studies  in 1918 in Warsaw, then continued studying mathematics and philosophy in Krak\'ow. In the school year 1925/26 he taught mathematics in the  state gymnasium in Rybnik. In 1926 he obtained PhD in philosophy at the Jagiellonian University on the basis of required exams and the thesis ``A priori knowledge in Bertrand Russell's epistemology", supervised by W\l adys\l aw Heinrich.  He translated into Polish some writings  of the dialectical materialistic philosopher Joseph Dietzgen   and poetic works by Schiller, Goethe and Heine.\\

Since 1916, Bilski was active in the Polish Socialist Party-Left. He established first communist structures  in Zgierz and organized worker demonstrations. In 1929, to avoid arrest by Polish authorities, he escaped to the USSR under the name of Stefan Biernacki. He supervised the Polish section of the Central Publishing House of the Nations of USSR in Moscow, then he lectured on the history of Polish labor movement in Kiev. He was arrested at the beginning of the Great Purge in 1934, sentenced as an alleged Polish spy  and executed  (\cite{Wie10}).\\

\subsection{Jan J\'ozef Le\'sniak (1901-1980)} 

\begin{wrapfigure}{l}{0.3\textwidth}

\begin{center}
\includegraphics[width=0.3\textwidth]{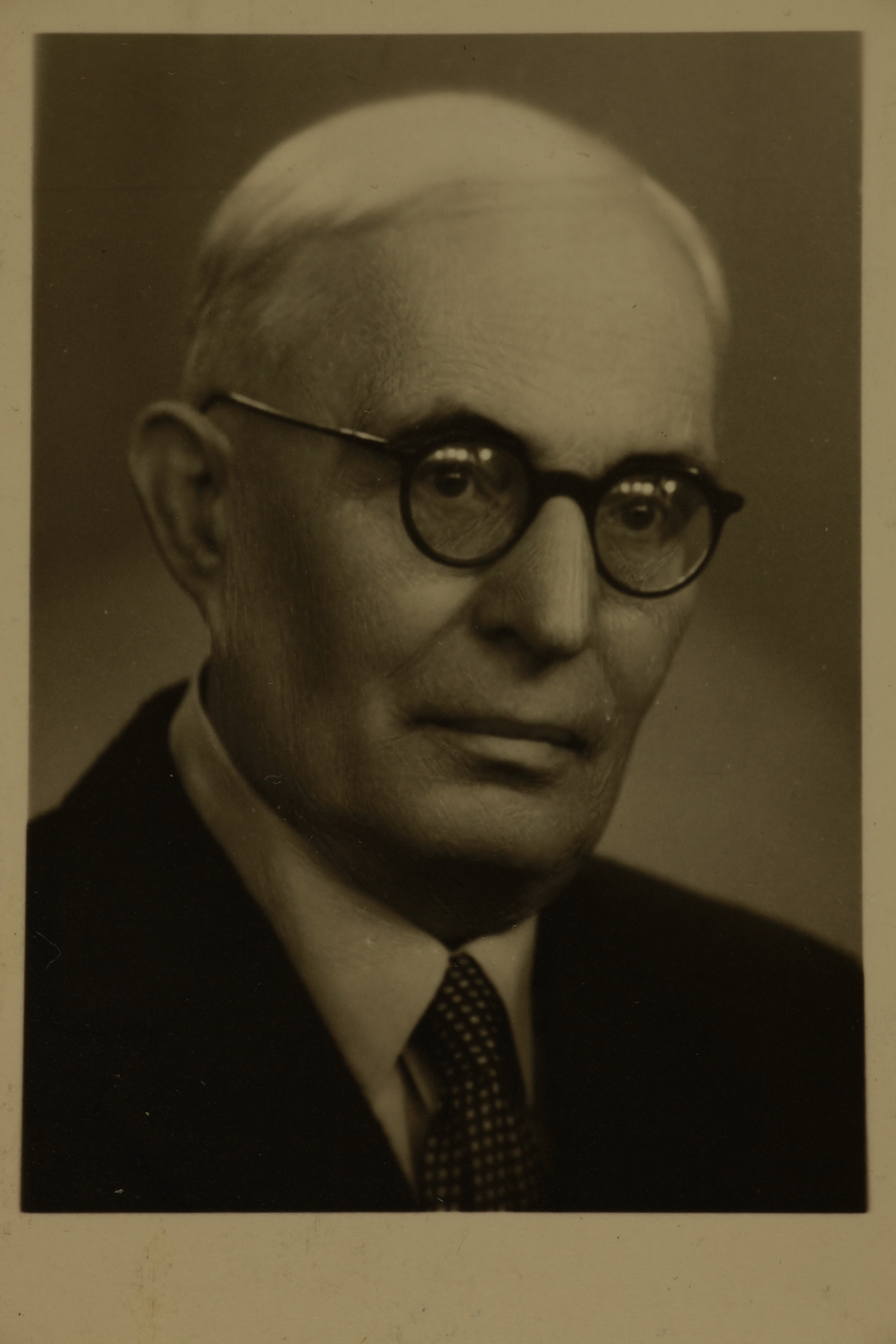}
\end{center}
%\caption{Jan J\'ozef Le\'sniak}
\end{wrapfigure}

Born in Ropczyce, he attended a gymnasium in Jas\l o. He passed his \emph{matura} with distinction in 1919, then he studied mathematics at the Jagiellonian University. In 1922 he became a scientific aide of the Mathematical Seminar, then he was made an assistant. In 1928 he passed the exam for high school
teaching licence and took a position in H. Sienkiewicz gymnasium in Krak\'ow. At the same time he gave contract lectures in issues of elementary mathematics at the Jagiellonian University and-- since 1930-- lectures on didactics of mathematics at the Pedagogical Institute of Studies. In 1939 he was arrested by Gestapo and taken to the concentration camps in Wi\'snicz and Auschwitz. Released in 1940, he returned to Krak\'ow, where he taught mathematics in the School for Commerce and Industry. After the war he re-assumed his position in a gymnasium and continued his lectures at the Jagiellonian University. In 1947 he obtained there his PhD degree on the basis of the thesis ``Methods of solving equations". In the same year he started working at the Pedagogical College (WSP) in Krak\'ow. In 1951 he submitted a habilitation thesis ``Educational values of instruction in mathematics and their fulfilment in high school" and became an extraodinary professor. In 1963 he was promoted to the level of a professor of mathematics at the Faculty of Mathematics, Physics and Chemistry at the Pedagogical College. In 1961, in the paper \cite{Le61},  he proposed a new formal approach to indefinite integrals.  His output counts over 20 items, including several books on elementary mathematics, theoretical arithmetics and functions of one variable.

\subsection{Stanis\l aw Go\l \c ab (1902-1980)}

\begin{wrapfigure}{L}{0.45\textwidth}

\begin{center}
\includegraphics[width=0.45\textwidth]{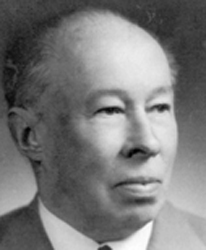}
\end{center}

\end{wrapfigure}

Stanis\l aw Go\l \c ab was born on July 26, 1902, in Travnik (Bosnia), where his father, Walenty, was a judge; his mother, Jadwiga ne\'e Skibi\'nska, was a teacher. In 1910 the family moved  to Krak\'ow. Go\l \c ab attended the V Gymnasium, where in 1920 he passed his maturity examination. One of his gymnasium teachers was Antoni Hoborski, later a professor of mathematics in the Academy of Mining, who also conducted lectures at the Jagiellonian University. In the years 1920-1924 Go\l \c ab  studied mathematics at the Jagiellonian University. In 1923 he completed the Pedagogical Study at the Jagiellonian University and in 1926 he passed examinations in mathematics and physics for candidates for high school teachers. In 1923 he started working as a deputy assistant to Hoborski, and after his graduation in 1924 was promoted to an assistant in the Chair of Mathematics at the  Academy of Mining.  Influenced by Hoborski (who wished to establish a scientific school in geometry in Krak\'ow), he chose geometry as his area of interest and started research in it. In the years 1925-1932 he published 14 papers. \\
In  1928, Go\l \c ab, with Hoborski's support, obtained a scholarship from the Division of Science of the Ministry of Religious Denominations and Public Education  for studying abroad. He went to Delft (the Netherlands) to work with J. A. Schouten. They wrote a paper together (\cite{SchGo30a}; Go\l \c ab also started working on his PhD thesis there. In 1929-1930 he got another scholarship, from the Fund for National Culture. He spent 3 months again in Delft, finishing his PhD thesis and writing the second part of his paper with Schouten (\cite{SchGo30b}). Then he went to Rome, where he learned absolute calculus from E. Bompiani and relativistic mechanics from T. Levi-Civita. The last three weeks he spent in Prague, in private discussions with L. Berwald and V. Hlavat\'y. His trip resulted in new publications (\cite{HlGo32}). In 1932, Go\l \c ab took part in the International Congress of Mathematicians in Z\"urich, giving two short talks there (\cite{Go32c}, \cite{Go32d}, \cite{ICM}). \\

Since 1930, Go\l \c ab conducted contract lectures at the Jagiellonian University. In 1931 he obtained the PhD degree at the Jagiellonian University on the basis of the thesis ``O uog\'olnionej geometrii rzutowej", which was presented in  Polish (and published as \cite{Go30a}) and was related to his joint work with Schouten. His supervisor was Stanis\l aw Zaremba; the referee was Witold Wilkosz. In 1932 he obtained habilitation at the Jagiellonian University on the basis of the paper ``Zagadnienia metryczne geometrii Minkowskiego", published in 1932 in the Proceedings of the Academy of Mining (\cite{Go32a}). The topic of his habilitation lecture was ``Metryka k\c atowa w og\'olnych przestrzeniach", chosen by the committee out of three topics proposed by the candidate (which he  developed in his congress talks and several other publications, including a joint paper with Adam Bielecki, \cite{BieGo45}). The examination questions concerned, among other things, Minkowski geometry, convex functions (Wilkosz), differential equations (Zaremba) and general metric spaces (Hoborski), and Go\l \c ab's answers were evaluated as very good.  The ministry extended Go\l \c ab's licence to lecture to the Faculty of Mining at the Academy of Mining, and he received the title of professor in 1939.\\

On November 6, 1939, Go\l \c ab was imprisoned along with other Krak\'ow professors and taken to the concentration camps, first in Sachsenhausen, then in Dachau. He tried to take care of his teacher Hoborski, but ultimately witnessed his death. Released \footnote{According to  \cite{Se03}, Wilhelm Blaschke claimed to intervene personally on Go\l \c ab's behalf. In private communication with the second author, Zofia Go\l \c ab-Meyer, Stanis\l aw Go\l \c ab's daughter, could not confirm whether this was the case. On the basis of documents in her possession, she brought to our attention the intervention of Hasso H\"arlen, one of Go\l \c ab's coauthors. This intervention is mentioned in \cite{Pie05}. } in December 1940 together with a group of junior scholars , he worked in forestry administration, taking part in clandestine teaching since 1943. After the liberation he continued working at the Academy of Mining, becoming an ordinary professor there in 1948.  When the State Mathematical Institute (later the Mathematical Institute of the Polish Academy of Sciences) was created, he became the head of the Division of Geometry there (he headed the division until his retirement). In the years 1950-55 he was a contract professor at the Pedagogical College  in Krak\'ow, and in 1954 he transferred from the Academy of Mining and Metalurgy to the Jagiellonian University (however, until 1962 he remained as a half-time employee at the Academy). He was the head of the Chair of Geometry and the dean of the Faculty of Mathematics, Physics and Chemistry of the Jagiellonian University. Because of his worthy behavior as the dean during the March Events of 1968 he was demoted from the position of the head of the Laboratory of Geometry by the communist authorities. He retired in 1972 and died in Krak\'ow on April 30, 1980. (See \cite{Kuch82}, \cite{Du12Go}, \cite{PB03Go}, \cite{GaPo00}.)\\

Go\l \c ab's  scientific output comprises over 250 publications, including two monographs: \textit{``Tensor Calculus"} (\cite{Go56}) and \textit{``Functional Equations in the Theory of Geometric Objects"} (joint with J. Acz\'el; \cite{GoAcz60}).  His main research interest was geometry, in particular clasical differential geometry, tensor calculus (under the influence of Schouten), spaces with linear or projective connection, Riemann, Minkowski and Finsler spaces, general metric spaces and other areas. He supervised 28 PhD degrees during his career. While working on the theory of geometric objects, in 1939 Go\l \c ab (\cite{Go39}) gave the first exact defintion of a pseudogroup of transformations (\cite{Kuch82}, \cite{Ho14}). He also promoted the methods of functional equations in geometry, starting with the so-called translation equation, for which he found a solution of class $C^1$ (\cite{Go38a}). But his interest in functional equations dates earlier. In 1937 (\cite{Go37}) he published a note in which he gave conditions characterizing the transformation $z \to \bar{z}$ of complex numbers as the only  solution to the multiplicative equation $f(uv)=f(u)f(v)$
(in the same journal issue, different conditions were given by A. Turowicz, \cite{Tu37}). He went on to obtain other results in  functional equations and supervise PhD candidates in this field, and is considered as the founding father of the Polish school of functional equations (\cite{Kucz73}).  In metric geometry, he proved e.g.  that, under certain mild assumptions, approximation of a point $x$ of  a Minkowski space by points in a given closed subset $C$  is unique if the distance function to $C$ is differentiable at $x$ (\cite{Go38b}; later on, many variants and refinements of this result were obtained by various authors, cf. \cite{Zaji83}). Another memorable result by Go\l \c ab in this field concerns the the unit disk in the plane with a  norm. It says that the perimeter of the unit disk can take any value between 6 and 8, and the extremal values are taken if and only if the disk is a regular hexagon or a parallelogram (\cite{Go32a}; \cite{AP05}).Go\l \c ab also worked on problems in applied mathematics, mainly in mining and geodesy (\cite{Bodz73}), in close collaboration with engineers. He was also interested in mathematics education and history of mathematics. He edited the monograph ``Studies in the history of the chairs in mathematics, physics and chemistry of the Jagiellonian University", published  in  celebration of the University's 600th anniversary  (\cite{Go64}). He  lectured with lucidity and liked working with students. In turn, they were fond of him as an educator and treated him as an authority.\\

\subsection{Zofia Krygowska (1904-1988)}  Anna Zofia Krygowska (ne\'e Czarkowska) was born in Lw\'ow. She grew up in Zakopane, where she finished gymnasium. In the years 1923-1927 she studied mathematics at the Jagiellonian University. Then she passed the teacher's licence examination and taught in elementary and high schools. In 1931 she obtained the master's degree from the Jagiellonian University. She became interested in issues concerning the process of teaching mathematics, in particular curriculum development and bringing new trends in mathematics into school-level education. She engaged in the activities of the Methodological Center for Mathematics, participated in many conferences and discussions and published articles on the subject of teaching mathematics (\cite{Cza36}). In 1937 she published a textbook ``Mathematics for the 1st grade of high school".\\

\begin{wrapfigure}{l}{0.45\textwidth}

\begin{center}
\includegraphics[width=0.45\textwidth]{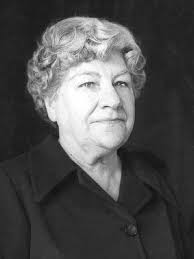}
\end{center}

\end{wrapfigure}
During the German occupation Krygowska was involved in clandestine teaching. Officially she worked as an accountant in a lumber company, while also travelling on behalf of the underground educational authorities  to organize, teach and coordinate illegal classes (especially in the Podhale region), risking her life. In 1945 she resumed her work as a high school teacher. In the years 1948-1951 she headed the Methodological Center. In 1950 she obtained the PhD degree at the Jagiellonian University on the basis of the thesis ``On the limits of rigor in the teaching of elementary geometry" prepared under the direction of Tadeusz Wa\.zewski. In this thesis she developed and studied an original system of axioms of geometry and proved its equivalence to Hilbert's system. (\cite{Turn05}). She also took a full-time position at the Pedagogical College in Krak\'ow. She began active efforts to modernize the education of teachers and to make didactics of mathematics a research discipline in its own right, while continuing her involvement in improving teaching of mathematics-- especially geometry-- in schools. In the years 67-71 she wrote or co-wrote five textbooks in geometry, which were subsequently adopted for high school use in Poland. She also published many articles and books on the issues of teaching of mathematics and its methodology.\\

In 1958 Krygowska became a professor in the newly created Unit of Methods of Teaching of Mathematics in the Pedagogical College in Krak\'ow. She initiated post-graduate doctoral programs in  didactics of mathematics. She supervised 22 PhD degrees and 4 habilitations in this field. She was a member of many national and international committees related to teaching mathematics, most importantly of  the Commission for the Study and Improvement of Mathematics Teaching (CIEAEM), founded in 1950. She served as a president, then as a honorary president, of this organization. She gave an address at the International Congress of Mathematicians in Nice in 1970. She was also instrumental in creating a series of lectures for teachers which were broadcast by Polish state television in the years 1968-77. She supported launching  a journal devoted to didactics of mathematics by the Polish Mathematical Society and  in 1982  the publication of  “Dydaktyka Matematyki” started.  She was an editor of this journal until her death. (\cite{No92}).\\

Krygowska became acquainted with modern  developments in set theory, topology, mathematical logic and other subjects as a student in Krak\'ow. She gave credit to Witold Wilkosz, who taught these subjects, for initiating her interest in them.  The new mathematical notions  as well as emphasis on rigor made their way to her school textbooks and were consistent with the trends (started in 1930s) of reformulating mathematics in a modern way, advocated by the Bourbaki group. Krygowska also acknowledged the role of psychology  in teaching. She was familiar with the work of Jean Piaget (whom she met at an international  conference in public education in Geneva in 1956) and followed some of his ideas in her paper ``Methodological and psychological basis for the activity-based method of teaching mathematics" (\cite{Kry57}).\\

Since her youth, Krygowska had a passion for mountaineering. She made a few first routes in the Tatra mountains, either in all-women teams or with her husband W\l adys\l aw Krygowski, a lawyer (\cite{RPP73}).   \\

\textbf{Acknowledgments:} We thank Dr. Zofia Pawlikowska-Bro\.zek, Dr. Zofia Go\l \c ab-Meyer, Professor Roman Duda  and Dr. Emelie Kenney for discussions on the topics related to Krak\'ow mathematics between the wars, as well as for their encouragement. The photos come from the archives of Z. Pawlikowska – Brożek and S. Domoradzki. They were obtained from private individuals or UJ Archives.  \\
The first author was partially supported by the Centre for Innovation and Transfer of Natural Sciences and Engineering Knowledge.


\begin{thebibliography}{9999999}

\bibitem{PANWaz} Archiwum PAN i PAU w Krakowie: Materials of Prof. T. Wa\.zewski

\bibitem{UJWaz} Archiwum UJ: T. Wa\.zewski's portfolios: personal, doctoral, nostrification, habilitation

\bibitem{Am12} L. Amoroso: Sopra un problema al contorno, \emph{Rendiconti del Circolo Matematico di Palermo}  33 (1): 75-85,(1912,) doi:10.1007/BF03015289, JFM 43.0453.03.



\bibitem{AP05} J. C. \'Alvarez Paiva,  A. Thompson: On the Perimeter and Area of the Unit Disc, \emph{The American Mathematical Monthly} Vol. 112, No. 2 (Feb., 2005), pp. 141-154 


\bibitem{BCMMS} I. Banic,  M. Crepnjak, M.  Merhar, U. Milutinovic, T. Sovic: Wa\.zewski's universal dendrite as an inverse limit with one set-valued bonding function, \emph{Glas. Mat. Ser. III} 48(68) (2013), no. 1, 137-165. MR3064249. ) Zbl 1276.54025

%\bibitem{Biel02} M. Bielawka: In the Shadow of the Master: Danuta Gierulanka, Phenomenology of Mathematics, \emph{Analecta Husserliana} 80:199-201 (2002) 

%\bibitem{Bie35} A Bielecki: \emph{O integralnem przedstawianiu $m$-wymiarowych powierzchni zawartych w $n$-wymiarowej przestrzeni euklidesowej zapomoc\c a funkcyj uwik\l anych} (PhD thesis), Dodatek do Rocznika Polskiego Towarzystwa Matematycznego, tom VII, Krak\'ow, 1935, 1-38 

%\bibitem{Bie48} A. Bielecki: 
 %Sur certaines conditions n\'ecessaires et suffisantes pour l'unicit\'e des solutions des syst\`emes d'\'equations %diff\'erentielles ordinaires et des \'equations au paratingent, \emph{Ann. Univ. Mariae Curie-Sk\l odowska} Sect. %A. 2 (1947), 49-106 (1948). MR0032885.  Zbl 0040.33401


%\bibitem{Bie56a} A. Bielecki: Une remarque sur la m\'ethode de Banach-Cacciopoli-Tikhonov dans la th\'eorie des \'equations diff\'erentielles ordinaires, \emph{Bull. Acad. Polon. Sci.}, Cl. III 4 (1956), 261-264. MR0082073.  Zbl 0070.08103
 


%\bibitem{Bie56b} A. Bielecki: Une remarque sur l'application de la m\'ethode de Banach-Cacciopoli-Tikhonov dans la th\'eorie de l'\'equation $s=f(x,y,z,p,q)$, \emph{Bull. Acad. Polon. Sci.}, Cl. III 4 (1956), 265-268. MR0082074. Zbl 0070.09004
  

%\bibitem{Bie56c} A. Bielecki: R\'eduction des axiomes de congruence de Hilbert, \emph{Bull. Acad. Polon. Sci.}, Cl. III 4 (1956), 321-324. MR0082106.  Zbl 0070.15802

%\bibitem{Bie57a} A. Bielecki: Extension de la m\'ethode du r\'etracte de T. Wa\.zewski aux \'equations au paratingent, \emph{Ann Univ. Mariae Curie-Sk\l odowska}, Sect. A 9 (1955), 37-61 (1957) MR0088620.  Zbl 0071.08902

%\bibitem{Bie57b} A. Bielecki: Sur l'ind\'ependance des axiomes d'incidence, d'ordre et de congruence de Hilbert, \emph{Ann Univ. Mariae Curie-Sk\l odowska}, Sect. A 9 (1955), 157-175 (1957) MR0082073.   Zbl 0078.12903

\bibitem{BieGo45} A. Bielecki, S. Go\l \c ab: Sur un probl\`eme de la m\'etrique angulaire dans les espaces de Finsler, \emph{Ann. Soc. Polon. Math.} 18 (1945), 134-144. MR0018965.  Zbl 0063.00380



%\bibitem{BieLe62} A. Bielecki, Z. Lewandowski: Sur certaines majorantes des fonctions holomorphes dans le cercle unit\'e, \emph{Colloq. Math.} 9 (1962), 299-303. MR0152650.  Zbl 0106.28303

%\bibitem{BMW39} A. Bielecki, M. Mathisson, J. W. Weyssenhoff: Sur un th\'eor\`eme concernant une transformation d'int\'egrales quadruples et int\'egrales curvilignes dans l'espace de Riemann, \emph{Bull. de l'Acad. Polon. des Sci. et des Lettr.}, S\'er.  A (1939), 22-28.  Zbl 0015.12002

%\bibitem{BW36} A. Bielecki, J. W. Weyssenhoff: Quaternions, $4$-dimensional rotations and Cayley's formula, \emph{Bull. de l'Acad. Polon. des Sci. et des Lettr.}, S\'er.  A (1936), 216-227.  Zbl 0015.12002


%\bibitem{BieZa36} A. Bielecki, S. K. Zaremba: Sur les points singuliers des syst\`emes de deux \'equations diff\'erentielles ordinaires, \emph{Ann. de la Soc. Polon. de Math.} 15 (1936), 135-139.  Zbl 0018.12302


\bibitem{Bodz73} J. Bodziony: Activity of Professor Stanis\l aw Go\l \c ab in the field of applied mathematics. Collection of articles dedicated to Stanis\l aw Go\l \c ab on his 70th birthday, I. \emph{Demonstratio Math.} 6 (1973), 45-49. MR0342350  Zbl 0274.01046

%\bibitem{BK14} J. Bourgain,  A. Kontorovich: On Zaremba's conjecture. \emph{Ann. of Math.} (2) 180 (2014), no. 1, 137-196. MR3194813



 %\bibitem{Ce05} A. Cellina:  A view on differential inclusions. \emph{Rend. Semin. Mat. Univ. Politec. Torino}  63 (2005), no. 3, 197-209. MR2201565

\bibitem{Ch98} J. J. Charatonik: History of continuum theory, in: E. Aull and R. Lowen (eds.), \emph{Handbook of the History of General Topology}, Volume 2, 703-786. 1998, Kluwer Academic Publishers.



%\bibitem{CP88} K. Ciesielski, Z. Pogoda: Conversation with Andrzej Turowicz, \emph{The Mathematical Intelligencer} 10.4 (1988), 13-20

\bibitem{CP12} K. Ciesielski, Z. Pogoda: On mathematics in Krak\'ow throughout centuries, \emph{Wiad. Mat.} 48(2), 2012, 313-324

%\bibitem{Co94} C. Corduneanu: Bielecki's method in the theory of integral equations, \emph{Ann. Univ. Mariae Curie-Sk\l odowska} Sect. A 39 (1994), 23-40

\bibitem{Cza36} Z. Czarkowska [later Krygowska]: O poj\c eciu granicy w nauczaniu matematyki w szkole [On the notion of limit in teaching mathematics at school], \emph{Matematyka i Szko\l a}, Warsaw, 1936 


%\bibitem{Da98} A. L. Dawidowicz, O pracach Zbigniewa \L omnickiego i Stanis\l awa Krystyna Zaremby ze statystyki matematycznej, in:  S. Fudali (ed.), \emph{Matematycy polskiego pochodzenia na obczy\'znie}, Wyd. Uniwersytetu Szczeci\'nskiego, Szczecin 1998, 201-209.

\bibitem{Dem14} S. S. Demidov: The origins of the Moscow school of the theory of functions, \emph{Technical Transactions}  1-NP, 2014, 73-84. 

\bibitem{Do13} S. Domoradzki: Cztery doktoraty z matematyki uzyskane przez Polak\'ow we Francji przed 1939 r., in: J. Be\v cva\v r, M. Be\v cva\v rova (eds.), 34. mezin\'arodni konference Historie Matematiky, Matfyzpress, Univerzita Karlova, Praha, 2013, 93-100

\bibitem{Do15} S. Domoradzki:  Mathematics in Polish territories in World War I, to appear

\bibitem{DS15} S. Domoradzki, M. Stawiska: Distinguished graduates in mathematics of Jagiellonian University in the interwar period. Part II:  1926-1939. This issue.

\bibitem{Du07} R. Duda: \emph{Lwowska Szko\l a Matematyczna}. Wydawnictwo Uniwersytetu Wroc\l awskiego, Wroc\l aw 2007.

\bibitem{Du12Go} R. Duda: Stanis\l aw Go\l \c ab, in: \emph{Matematycy XIX i XX wieku zwi\c azani z Polsk\c a}, Wyd. Uniwersytetu Wroc\l awskiego, Wroc\l aw 2012, 151-153

%\bibitem{Du12Tu} R. Duda: Andrzej Turowicz, in: \emph{Matematycy XIX i XX wieku zwi\c azani z Polsk\c a}, Wyd. Uniwersytetu Wroc\l awskiego, Wroc\l aw 2012, 486-488

%\bibitem{Du12Za} R. Duda: Stanis\l aw Krystyn Zaremba, in: \emph{Matematycy XIX i XX wieku zwi\c azani z Polsk\c a}, Wyd. Uniwersytetu Wroc\l awskiego, Wroc\l aw 2012, 524-525

\bibitem{Dy00} J. Dybiec: Uniwersytet Jagiello\'nski 1918-1939. Nak\l adem Polskiej Akademii Umiej\c etno\'sci, Krak\'ow 2000, 759 pp.

\bibitem{GaPo00} J. Gancarzewicz, Z. Pogoda: \emph{Stanis\l aw Go\l \c ab (1902-1980).} In: B. Szafirski (ed.): Z\l ota Ksi\c ega wydzia\l u matematyki i Fizyki UJ, Krak\'ow 2000, p. 357-362

%\bibitem{Gie58} D. Gierulanka: O przyswajaniu sobie poj\c e\'c geometrycznych [On Acquiring Geometrical Notions], Pa\'nstwowe Wydawn. Naukowe, Warsaw, 1958

%\bibitem{Gie62} D. Gierulanka: Zagadnienie swoisto\'sci poznania matematycznego [The Problem of Specificity of the Mathematical Cognition], Pa\'nstwowe Wydawn. Naukowe, Warsaw, 1962

%\bibitem{Fe84} Fedorchuk, V. V.: Some geometric properties of covariant functors, \emph{Uspekhi Mat. Nauk} 39 (1984), no. 5(239), 169–208 MR0764014 (86j:54002)


\bibitem{Go30a} S. Go\l \c ab: \"Uber verallgemeinerte projektive Geometrie.  \emph{Prace mat.-fiz.} 37, 91-153 (1930). JFM 56.1194.01

\bibitem{Go30b} S. Go\l \c ab: Sopra le connessioni lineari generali. Estensione d'un teorema di Bompiani nel caso pi\'u generale. \emph{Annali di Mat.} (4) 8, 283-291 (1930) JFM 56.0629.01

\bibitem{Go32a} S. Go\l \c ab: Quelques probl\`emes m\'etriques de la g\'eometrie de Minkowski, \emph{Trav. de l'Acad. Mines Cracovie} 6 (1932), pp. 179

\bibitem{Go32c} S. Go\l \c ab: \"Uber die M\"oglichkeit einer absoluten Auszeichnung der Gruppe von Koordinatensystemen in verschieden R\"aumen,  internationalen Mathematikerkongress, Z\"urich, 1932

\bibitem{Go32d} S. Go\l \c ab: Einige Bemerkungen \"uber Winkelmetrik in Finslerschen R\"aumen, internationalen Mathematikerkongress, Z\"urich, 1932

\bibitem{Go37} S. Go\l \c ab: Sur une d\'efinition axiomatique des nombres conjugues pour les nombres complexes ordinaires, \emph{Opuscula Math.} 1, 1-11 (1937). Zbl 0019.01002

\bibitem{Go38a} S. Go\l \c ab: \"Uber eine Funktionalgleichung der Theorie der geometrischen Objekte, \emph{Wiadom. Mat.} 45 (1938), 97-137

\bibitem{Go38b} S. Go\l \c ab: Sur la fonction repr\'esentant la distance d'un point variable \`a un ensemble fixe, \emph{C. R. Acad. Sci., Paris} 206, 406-408 (1938).  Zbl 0018.17205

\bibitem{Go39} S. Go\l \c ab: 
\"Uber den Begriff der ``Pseudogruppe von Transformationen".  
\emph{Math. Ann.}, Berlin, 116, 768-780 (1939). JFM 65.1122.03

\bibitem{Go56} S. Go\l \c ab: \emph{Rachunek Tensorowy}. Bibl. Mat. 11, Warszawa, 1956

\bibitem{Go64} S. Go\l \c ab: Zarys dziej\'ow matematyki w UJ w XX w.
W: \emph{Studia z dziej\'ow katedr Wydzia\l u Matematyki, Fizyki i Chemii}, Krak\'ow, UJ 1964, 75-86


\bibitem{GoAcz60} S. Go\l \c ab, J. Acz\'el: \emph{Funktionalgleichungen der Theorie der geeometrischen Objekte}. Monografie Mat. 39, Warszawa, 1960  

%\bibitem{GoTu70} H. G\'orecki, A. Turowicz: \emph{Optimal control. A survey on mathematical methods. %(Sterowanie optymalne. Przegl\c ad metod matematycznych.)}. (Polish) 
%Monografie zagadnie\'n elektrotechniki teoretycznej. Warszawa-Wroc\l aw: Pa\'nstwowe Wydawnictwo Naukowe. %(1970).Zbl 0235.49002

\bibitem{Gor11} L. G\'orniewicz:  Differential inclusions--the theory initiated by Cracow Mathematical School. \emph{Ann. Math. Sil.} No. 25 (2011), 7-25 (2012) MR2986421

\bibitem{HlGo32} V. Hlavat\'y, S. Go\l \c ab: Zur Theorie der Vektor- und Punktkonnexion.  \emph{Prace Mat.-Fiz.} 39, 119-130 (1932). JFM 58.0762.01

\bibitem{Ho14} C. Hollings: \emph{Mathematics across the Iron Curtain.
A history of the algebraic theory of semigroups.} History of Mathematics, 41. American Mathematical Society, Providence, RI, 2014. xii+441 pp. ISBN: 978-1-4704-1493-1 MR3222721 Zbl 06329297

\bibitem{HLP14} E. Hrushovski, F. Loeser, B. Poonen: Berkovich spaces embed in Euclidean spaces, \emph{ L'Enseignement Math\'ematique} 60, 273-292 (2014) 


\bibitem{ICM} ICM Proceedings since 1893, International Mathematical Union, www.mathunion.org/ICM/

\bibitem{KaNi28} S. Kaczmarz, W. Nikliborc: Sur les suites de fonctions convergentes en moyenne, \emph{Fundamenta Mathematicae} 11 (1928), 151-168

%\bibitem{KaTu} S. Kaczmarz, A. Turowicz: Sur l'irrationalit\'e des int\'egrales ind\'efinies, \emph{Studia Math.}, Lw\'ow, 8, 129-134 (1939).JFM 65.0212.02

%\bibitem{KKZ04} W. Kaczor, T. Kuczumow, W. Zygmunt: Adam Bielecki (1910-2003), \emph{Roczniki Pol. Tow. Mat. Seria II: Wiadomo\'sci Matematyczne} XL (2004), 213-228. MR2337182


%\bibitem{Koj08a} M. Kojdecki: Instytut Matematyki i Kryptologii. Zarys historii 1968-2008.
%strona.wcy.wat.edu.pl/index.php?option=com_content\&task=view\&id=65\&Itemid=59
%Na podstawie ksi\c a\.zki \emph{Wydzia\l\ Cybernetyki Wojskowej Akademii Technicznej 1968-2008.}
%Praca zbiorowa pod redakcj\c a naukow\c a Wojciecha W\l odarkiewicza. Warszawa 2008 ISBN 978-83-89399-96-0

%\bibitem{Koj08b} M. Kojdecki: Roman Leitner. \emph{G\l os Akademicki WAT}, June 2008, no. 6(147), p.11 %www.promocja.wat.edu.pl/Glos_Akademicki/Glos_PDF/2008/g



 %\bibitem{Kon13} A. Kontorovich: From Apollonius to Zaremba: local-global phenomena in thin orbits. \emph{Bull. Amer. Math. Soc. (N.S.)} 50 (2013), no. 2, 187-228 MR3020826

\bibitem{KrUJ2122} Kronika Uniwersytetu Jagiellon\'nskiego za lata szkolne 1921/22 i 1922/23. 

\bibitem{Kry57} Z. Krygowska: Metodologiczne i psychologiczne podstawy czynno\'sciowej metody nauczania matematyki [Methodological and psychological foundations of the activity-based method in teaching mathematics], in: \emph{The Tenth Anniversary of the Pedagogical College in Krak\'ow 1946-1956}, PWN Krak\'ow, 1957 

\bibitem{KsSz3541} Copies of the manuscript and  the typescript of the Scottish Book (in Polish and English),  Stefan Banach's wortal, kielich.amu.edu.pl, accessed August 09, 2015

\bibitem{Kuch73} M. Kucharzewski:
Scientific achievements of Professor Stanislaw Go\l \c ab in the domain of geometry. Collection of articles dedicated to Stanisl\ aw Go\l \c ab on his 70th birthday, I. \emph{Demonstratio Math.} 6 (1973), 19-38.  MR0342348 Zbl 0274.01044

\bibitem{Kuch82}   M. Kucharzewski:  Stanislaw Go\l \c ab-life and work. \emph{Aequationes Math.} 24 (1982), no. 1, 1-18. MR0698112  Zbl 0504.01017


\bibitem{Kucz73} M. Kuczma: Activity of Professor Stanis\l aw Go\l \c ab in the theory of functional equations. Collection of articles dedicated to Stanis\l aw Go\l \c ab on his 70th birthday, I. \emph{Demonstratio Math.} 6 (1973), 39-44. MR0342349  Zbl 0274.01045

%\bibitem{Ku02} L. Kusak: Danuta Gierulanka - \.zycie i tw\'orczo\'s\'c, \emph{Ruch Filozoficzny} 3 (3) (2002)

  
\bibitem{Le61}  J. Le\'sniak, On the definitions of indefinite integrals, \emph{Wy\.z. Szko\l . Ped. Krak\'ow. Rocznik Nauk.-Dydakt.} no. 13, Prace Mat. 1961, 137-152, MR0308339 (46 \#7453)

 
%\bibitem{LCh95} G. \L oskiewicz, B. Choczewski: Life and work of Tadeusz Rachwa\l\ (1914-1992) %\emph{Opuscula Math.} No. 15 (1995), 119-125. MR1365939 (96h:01035)

%\bibitem{MS00} J. Madey, M. Sys\l o: Pocz\c atki informatyki w Polsce. \emph{Informatyka} 9,10 (2000).
%English version: IEEE Annals of History of Computing. 
%Regnecentralen in Eastern Europe. 
%http://datamuseum.dk/site_dk/rc/NIB/kap14.shtml

%\bibitem{Ma71} Maturzy\'sci 71 (collective):  ``Ciocia" i inni.  [``Auntie" and others.] In: \emph{W Krakowie %przy Studenckiej. Ksi\c ega Pami\c atkowa wydana z okazji Jubileuszu 130-lecia V Liceum Og\'olnokszta\l c\c %acego im. A. Witkowskiego w Krakowie 1871-2001}, ed. Grzegorz Ma\l achowski, Oficyna Wydawnicza Text, %Krak\'ow 2001


\bibitem{MO87} G. H.Meisters, C. Olech: A poly-flow formulation of the Jacobian conjecture, \emph{Bull. Polish Acad. Sci. Math.} 35 (1987), no. 11-12, 725–731. MR0961711 (89j:13005)

\bibitem{Ni25} L. Nikliborc: Sur les fonctions hyperharmoniques,  
\emph{C. R. Acad. Sci. Paris }180, 1008-1011 (1925). JFM 51.0364.02

\bibitem{Ni26} L. Nikliborc: Sur les fonctions hyperharmoniques,   
{C. R. Acad. Sci. Paris} 182, 110-112 (1926). JFM 52.0498.02

\bibitem{No92} B. J. Nowecki: Anna Zofia Krygowska, \emph{Educational Studies in Mathematics}, Vol. 23, No. 2 (Apr., 1992), pp. 123-137


%\bibitem{NT55}  J.Nowi\'nski, S.  Turski: Solution of the Lam\'e problem for a heterogeneous cylinder by means of the ARR differential equation analyser. (Polish) Arch. Mech. Stos. 7 (1955), 419-424. MR0074254 

\bibitem{Ol76} C. Olech: The achievements of Tadeusz Wa\.zewski in the mathematical theory of optimal control. (Polish) \emph{Wiadom. Mat.} (2) 20 (1976), no. 1, 66–69. MR0433279 (55 \# 6257)

%\bibitem{Ol10} C.  Olech: Selections from the curriculum vitae: happy thirteens. (Polish) Wiad. Mat. 46 (2010), no. 2, 215-224. MR2905712

\bibitem{Pa03} T. Pankiewicz: \emph{Apteka w getcie krakowskim}, Wydawnictwo Literackie, Krak\'ow, 2003, str.46-48. 

\bibitem{PB03Go} Z. Pawlikowska-Bro\.zek:  Stanis\l aw Go\l \c ab, in: \emph{S\l ownik Biograficzny Matematyk\'ow Polskich}, S. Domoradzki, Z. Pawlikowska-Bro\.zek, D. W\c eglowska (Eds.), Wyd. PWSZ, Tarnobrzeg 2003, 62-63

%\bibitem{PB03Tu} Z. Pawlikowska-Bro\.zek:  Andrzej Turowicz, in: \emph{S\l ownik Biograficzny Matematyk\'ow Polskich}, S. Domoradzki, Z. Pawlikowska-Bro\.zek, D. W\c eglowska (Eds.), Wyd. PWSZ, Tarnobrzeg 2003, 246-247.

%\bibitem{PB03Za} Z. Pawlikowska-Bro\.zek:  Stanis\l aw K. Zaremba, in: \emph{S\l ownik Biograficzny Matematyk\'ow Polskich}, S. Domoradzki, Z. Pawlikowska-Bro\.zek, D. W\c eglowska (Eds.), Wyd. PWSZ, Tarnobrzeg 2003, 268-269.

%\bibitem{Pe00a}  A. Pelczar: O matematyce i matematykach  w  Uniwersytecie Jagiello\'nskim,  in:  B. Szafirski (ed.), Z\l ota Ksi\c ega UJ, Wydzia\l\ Matematyki i Fizyki, 214-237. 

\bibitem{Pe00b} A.Pelczar: Tadeusz Wa\.zewski (1896-1972) uczony i nauczyciel, [in:] Wydzia\l\
Matematyki i Fizyki, Z\l ota Ksi\c ega, Uniwersytet Jagiello\'nski, 600-lecie odnowienia Akademii
Krakowskiej, Krak\'ow 2000, 341-356.

\bibitem{Pie05} H. Pierzcha\l a: Pomocne d\l onie Europejczyk\'ow (1939-1944), Wydawnictwo i Poligrafia Pijar\'ow, Krak\'ow 2005.
	

%\bibitem{PG78}  A. Pelczar, H. G\'orecki: O dzia\l alno\'sci naukowej Profesora Andrzeja Turowicza, \emph{Wiadom. Mat.} 21 (1978), 15-24

%\bibitem{PGM05} A. Pelczar, H. G\'orecki, W. Mitkowski:  The centenary of the birth of Professor A. B. Turowicz. (Polish), \emph{Wiadom. Mat.} 41 (2005), 151-164. MR2340251

\bibitem{RPP73} Z. Radwa\'nska-Paryska, W. H. Paryski: Encyklopedia Tatrza\'nska, Warszawa 1973.

%\bibitem{Ra73} T. Rachwa\l : \emph{Geometria wykre\'slna}

%\bibitem{Ra01} T. Rachwa\l : Jak zosta\l em matematykiem. In: \emph{W Krakowie przy Studenckiej. Ksi\c ega Pami\c atkowa wydana z okazji Jubileuszu 130-lecia V Liceum Og\'olnokszta\l c\c acego im. A. Witkowskiego w Krakowie 1871-2001}, ed. Grzegorz Ma\l achowski, Oficyna Wydawnicza Text, Krak\'ow 2001




%\bibitem{RT35} A. Rosenblatt, S. Turski: Sur les coefficients des s\'eries de puissances univalentes dans le cercle unite,   \emph{C. R. Acad. Sci., Paris} 200, 1270-1272 (1935).Zbl 0011.26101

%\bibitem{RT36b} A. Rosenblatt, S. Turski: Sur la repr\'esentation conforme de domaines plans,   \emph{Bull. Sci. Math.}, II. Ser. 60, 309-320 (1936).Zbl 0015.11602

%\bibitem{RT36c} A. Rosenblatt, S. Turski: Sur la repr\'esentation conforme de domaines plans,   \emph{C. R. Acad. Sci., Paris} 202, 899-901 (1936).Zbl 0013.31203


\bibitem{Ry12}  K. P. Rybakowski:  A note on Wa\.zewski principle and Conley index. \emph{Wiad. Mat.} 48 (2012), no. 2, 223-237. MR2986196

\bibitem{SchGo30a} J. A.Schouten, S. Go\l \c ab:
\"Uber projektive \"Ubertragungen und Ableitungen.  
\emph{M. Z.} 32, 192-214 (1930). JFM 56.0627.02

\bibitem{SchGo30b}  J. A.Schouten, S. Go\l \c ab: \"Uber projektive \"Ubertragungen und Ableitungen. II.  
\emph{Annali di Mat.} (4) 8, 141-157 (1930). JFM 56.0628.01

\bibitem{Se03} S. L. Segal: \emph{Mathematicians under the Nazis}. Princeton University Press, Princeton, NJ, 2003. xxiv+530 pp. ISBN: 0-691-00451-X MR1991149  Zbl 1028.01006


 \bibitem{Srz11} R. Srzednicki: Metoda retraktowa Wa\.zewskiego, in: Tadeusz Wa\.zewski 1896-1972, PAU, Archiwum Nauki PAN i PAU, \emph{W S\l u\.zbie Nauki}  17, 2011, s. 33-36

\bibitem{Sl48} W. \'Slebodzi\'nski: W\l adys\l aw Nikliborc et son oeuvre scientifique. (French) \emph{Colloquium Math.} 1, (1948). 322-330. MR0029836

%\bibitem{Sre07} B. \'Sredniawa: Autobiografia. Wspomnienie fizyka teoretyka w 90-lecie urodzin i w 60-lecie promocji doktorskiej. \emph{Kwartalnik Historii Nauki i Techniki} 52/3-4 (2007), 7-19 

\bibitem{Turn05} S. Turnau: Profesor Zofia Krygowska jako matematyk, \emph{Roczniki Polskiego Tow. Matematycznego Seria V: Dydaktyka Matematyki}, 28 (2005), 79-82


\bibitem{TuAU} A. Turowicz: tape recordings; made available courtesy of Dr. Z. Pawlikowska-Bro\.zek

\bibitem{Tu37} A. Turowicz: Sur une d\'efinition axiomatique des nombres conjugu\'es. \emph{Opuscula Math. } 1, 13-16 (1937)

%\bibitem{Tu48} A. Turowicz: Sur les fonctionnelles continues et multiplicatives,   \emph{Ann. Soc. Polon. Math.} 20, 135-156 (1948). Zbl 0031.22001

%\bibitem{Tu67} A. Turowicz: \emph{Geometria zer wielomian\'ow.} (Polish) [Geometry of zeros of polynomials], Pa\'nstwowe Wydawnictwo Naukowe, Warsaw 1967, 155 pp. MR0229806 (37 \#5372) 

%\bibitem{Tu85} A. Turowicz: \emph{Teoria macierzy, Wykłady na Studium Doktoranckim w zakresie automatyki i elektroniki w roku 1970/71},  notes by  Wojciech Mitkowski, ed. IV, AGH, Kraków 1985 (ss. 240). 

%\bibitem{Tu95}  A. Turowicz:  On a proof of the Weierstrass-Stone theorem. (Polish), \emph{Wiadom. Mat.} 31 (1995), 149–150. MR1392944 (97c:01031). Zbl 0872.01015





\bibitem{Vi23} L. Vietoris, Kontinua zweiter Ordnung, \emph{Monatsh. Math. und Phys.} 33 (1923), 49-62

\bibitem{Wa23a} T. Wa\.zewski: Sur un continu singulier, \emph{Fund. Math.} 4 (1923a), 214-245

\bibitem{Wa23b} T. Wa\.zewski: Sur les courbes de Jordan ne renfermant aucune courbe simple ferm\'ee de Jordan,
\emph{Annales de la Soci\'et\'e Polonaise de Math\'ematique} 2(1923), 49-170.

\bibitem{Wa27} T. Wa\.zewski: Kontinua prostowalne w zwi\c azku z funkcjami i odwzorowaniami absolutnie ci\c ag\l ymi, \emph{Dodatek do Roczn. Pol. Tow. Mat.} (1927), 9-49.

\bibitem{Wa31} T. Wa\.zewski: Remarque sur un th\'eor\`eme de M. Bielecki, \emph{Ann. Soc. Polon. Math.}, 10(1931), 42-44

\bibitem{Wa45} T. Wa\.zewski: Sur une m\'ethode approximative de M. Rappaport concernant la trisection d'un angle, \emph{Ann. Soc. Polon. Math.}, 18(1945), 164 (Comptes Rendus des s\'eances des sections, sect.
Cracovie, 29.V.1945)

\bibitem{Wa47}  T. Wa\.zewski: Sur un principe topologique de l'examen de l'allure asymptotique des int\'egrales des \'equations diff\'erentielles ordinaires, \emph{Ann. Soc. Polon. Math.}, 20(1947), 279-313

\bibitem{Wa48} T. Wa\.zewski: Sur l'\'evaluation du domaine d'existence des fonctions implicites r\'eelles ou complexes, \emph{Ann. Soc. Polon. Math.} 20 (1947), 81–120 (1948). MR0026100 (10,106h)

\bibitem{Wa48l}  T. Wa\.zewski: Elementy rachunku r\'o\.zniczkowego i ca\l kowego, Krak\'ow 1948, 1-127 (lecture notes),  wyd. Bratnia Pomoc Student\'ow

\bibitem{Wa49} T. Wa\.zewski; Quelques d\'emonstrations uniformes pour tous les cas du th\'eor\`eme de l'H\^ospital. G\'en\'eralisations, \emph{Prace Mat. Fiz.}, 47, 117-128.

%\bibitem{Wa54}  T. Wa\.zewski: Une g\'en\'eralisation des t\'eor\`emes sur les accroissements finis au cas des espaces de Banach et application \`a la g\'enéralisation du th\'eor\`eme de l'H\^opital.
%Ann. Soc. Polon. Math. 24 (1951), no. 2, 132–147 (1954).MR0060728 (15,717g)

\bibitem{Wa61} T. Wa\.zewski: Sur la s\'emicontinuit\'e inf\'erieure du ``tendeur'" d'un ensemble compact variant d'une fa\c con continue, \emph{Bull. Acad. Polon. Sci. S\'er. Sci. Math. Astronom. Phys.} 9 (1961) 869-872. MR0132808 (24 \# A2644)




%\bibitem{Wa35} MR1556907
%Wa\.zewski, T.;
%Sur les matrices dont les éléments sont des fonctions continues.
%Compositio Math. 2 (1935), 63–68.

\bibitem{WaZa35}  T. Wa\.zewski,   S.K.Zaremba: Les ensembles limites des integrales des syst\`emes
d'\'equations diff\'erentielles, \emph{Ann. Soc. Polon. Math.}, 14(1935), 181 {Comptes Rendus des
s\'eances des sections, sect. Cracovie, 2.XII,1935}

\bibitem{WaZa36} T. Wa\.zewski,  S.K.Zaremba: Les ensembles de condensation des caract\'eristiques d'un
syst\`eme d' \'equations diff\'erentielles ordinaries, \emph{Ann. Soc. Polon. Math.}, 15(1936), 24-33

\bibitem{Wie10} M. Wierzbowski: \emph{Ilustrowany Tygodnik Zgierski} 25 (22-28.06.2010), cms.miasto.zgierz.pl/index.php?page=od-nazwisk-lub-imion

%\bibitem{Wit} E. Wittbrod: Zarys historii Politechniki w Gda\'nsku po 1945 r., pg.edu.pl/uczelnia/historia/po-1945


\bibitem{Wo39} M. Wojdys\l awski, Retractes absolus et hyperspaces des continu, 
\emph{Fund. Math.} 32 (1939), 23-25

%\bibitem{Zaj12} M. A.  Zaj\c ac: \emph{Fenomenologiczne w\c atki w psychologii poznania. Badania Danuty Gierulanki nad przyswajaniem poj\c e\'c i rozumieniem tekstu.}Wydawnictwo Uniwersytetu \'Sl\c askiego, Katowice 2012

\bibitem{Zaji83} L. Zaji\v cek, Differentiability of the distance function and points of multivaluedness of the metric projection in Banach space, \emph{Czechoslovak Math. J.}  33(108) (1983), no. 2, 292-308.


%\bibitem{Za35} S. K. Zaremba: O r\'ownaniach paratyngensowych [On paratingent equations], Dodatek do rocznika PTM, 9(1935), rozprawa doktorska (PhD thesis),

%\bibitem{Za36} S. K. Zaremba: Sur les \'equations au paratingent, \emph{Bull. Sci. Math.} 60 (2) (1936), 139-160.

%\bibitem{Za39} S. K. Zaremba: Sur l'allure des int\'egrales d'une \'equation diff\'erentielle  ordinaire du premier ordre dans le voisinage de l'int\'egrale singuli\`ere,  \emph{Bulletin International de l'Academie des Sciences de Cracovie}, ser. A, 1939

%\bibitem{Za71} S. K. Zaremba:  La m\'ethode des ``bons treillis" pour le calcul des int\'egrales multiples. In: \emph{Proc. Sympos., Univ. Montreal}, Montreal, Que., 1971, pp. 39-119. Academic Press, New York, 1972.


%\bibitem{ZaXX} S. K. Zaremba, personally written CV, copy courtesy of dr Z. Pawlikowska-Bro\.zek,  Krak\'ow. 


\end{thebibliography}
\end{document}